\documentclass[11pt,letterpaper]{article}
\usepackage{amsthm}
\usepackage{amsmath}
\usepackage{amsfonts}
\usepackage{latexsym}
\usepackage{geometry}
\usepackage{enumerate}
\usepackage{lipsum}

\emergencystretch=\maxdimen
\title{On equivalence between noncollapsing and bounded entropy for ancient solutions to the Ricci flow}
\author{Yongjia Zhang}
\numberwithin{equation}{section}
\begin{document}
\maketitle

In \cite{zhang2017entropy} we showed that for an ancient solution to the Ricci flow with nonnegative curvature operator, assuming bounded geometry on one time slice, bounded entropy implies noncollapsing on all scales. In this paper we prove the implication in the other direction, that for an ancient solution with bounded nonnegative curvature operator, noncollapsing implies bounded entropy. Hence we prove Perelman's assertion under the assumption of bounded geometry on one time slice. In particular, for ancient solutions of dimension three, we need only to assume bounded curvature. We also establish an equality between the asymptotic entropy and the asymptotic reduce volume, which is a result similar to Xu \cite{xuequation}, where he assumes noncollapsing and the Type I curvature bound.

\tableofcontents

\section{Introduction}
\newtheorem{Main_Theorem}{Theorem}[section]
\newtheorem{Main_Theorem_2}[Main_Theorem]{Theorem}
\newtheorem{Main_Corollary}[Main_Theorem]{Theorem}
\newtheorem{Gaussian_bound}[Main_Theorem]{Corollary}

In Perelman \cite{perelman2002entropy}, he asserted that for ancient solutions to the Ricci flow with nonnegative curvature operator, bounded entropy is equivalent to noncollapsing on all scales. He wrote
\begin{quotation}
...require that $g_{ij}(t)$ to be $\kappa$-noncollapsed on all scales... It is not hard to show that this requirement is equivalent to a uniform bound on the entropy $S$...
\end{quotation}
Here he was referring to the \textit{pointed entropy} which we will define in the next section. Though he claimed this equivalence to be obvious, we are not aware of a complete proof yet. In this paper, we shall prove that on an ancient noncollapsed solution to the Ricci flow with bounded nonnegative curvature operator, the pointed entropy is uniformly bounded. Combined with our previous work \cite{zhang2017entropy}, we prove Perelman's assertion with an additional assumption of bounded geometry on one time slice, that is, an upper bound of the curvature and a positive lower bound of the volume of all unit balls.
\\

\begin{Main_Theorem} \label{Main_Theorem}
Let $(M^n,g(\tau))_{\tau\in[0,\infty)}$ be an ancient solution to the backward Ricci flow with nonnegative curvature operator, assume that $\displaystyle\sup_{x\in M}R(x,0)<\infty$ and $\displaystyle \inf_{x\in M}\operatorname{Vol}(B_{g(0)}(x,1))>0$. Then $(M,g(\tau))_{\tau\in[0,\infty)}$ is noncollapsed on all scales if and only if there exists $x_0\in M$, such that $\bar{W}(x_0)>-\infty$, where $\bar{W}(x_0)$ is the asymptotic entropy defined in (\ref{eq:def_W}).
\end{Main_Theorem}

We remark here that by Hamilton's trace Harnack \cite{hamilton1993harnack}, the condition $\sup_{x\in M}R(x,0)<\infty$ is equivalent to $\sup_{(x,\tau)\in M\times[0,\infty)}R(x,\tau)<\infty$. In the above statement, the noncollapsing assumption refers to either weak noncollapsing or strong noncollapsing, they are equivalent in our case as remarked by Perelman \cite{perelman2002entropy}. A backward Ricci flow $(M^n,g(\tau))$ is called \textit{weakly $\kappa$-noncollapsed} if for any $(p,\tau)$ and $r>0$, $|Rm|\leq r^{-2}$ on $B_{g(\tau)}(p,r)\times[\tau,\tau+r^2]$ implies $\operatorname{Vol}(B_{g(\tau)}(p,r))\geq\kappa r^n$; is called \textit{strongly $\kappa$-noncollapsed} if for any $(p,\tau)$ and $r>0$, $R\leq r^{-2}$ on $B_{g(\tau)}(p,r)$ implies $\operatorname{Vol}(B_{g(\tau)}(p,r))\geq\kappa r^n$. In general cases, these notions are not equivalent. One may use a compact Ricci flat manifold whose Riemann curvature vanishes nowhere to construct a static ancient solution. It is easy to see such an ancient flow is weakly noncollapsed but not strongly noncollapsed. Moreover, according to the computation in our previous work, on such an ancient flow the entropy must be unbounded; see Lemma 3.1 of \cite{zhang2017entropy}. So if one aims at proving the general statement, the strong noncollapsing should be the correct notion. Notice that in this paper, by noncollapsing we always mean that there exists a $\kappa>0$, such that the $\kappa$-noncollapsing assumption holds on all scales.
\\

It turns out that on a noncollapsed ancient solution with bounded nonnegative curvature operator, the asymptotic entropy of the ancient solution coincides with the entropy of the asymptotic Ricci shrinker. Combining this with a result of Carrillo and Ni \cite{carrillo2009sharp}, we obtain an equality between the asymptotic entropy and the \textit{asymptotic reduced volume}, which is the same conclusion as Xu \cite{xuequation} with his Type I and noncollapsing condition replaced by the nonnegative curvature operator and the bounded geometry assumption.
\\

\begin{Main_Corollary} \label{Main_Corollary}
Under the same assumption as Theorem \ref{Main_Theorem}, we have that
\begin{eqnarray} \label{eq:equivalence}
\bar{W}(x_0)=\log\bar{V}(x_0),
\end{eqnarray}
where $\bar{V}(x_0)$ is the asymptotic reduced volume defined in (\ref{def:asymp_rvol}). Notice if the ancient flow is collapsed, then we have both sides of (\ref{eq:equivalence}) are $-\infty$.
\end{Main_Corollary}

\bigskip

Now we have that under the noncollapsing assumption, (\ref{eq:equivalence}) holds in two cases, either in the Type I case or in the bounded nonnegative curvature operator case, which are the only cases we know the existence of the asymptotic Ricci shrinker (for the asymptotic Ricci shrinker in Type I case, see Naber \cite{naber2010noncompact}). As long as the asymptotic Ricci shrinker exists, the entropy and the reduced volume descend at the same time to it, on which these two quantities coincide in some certain way; see Carrillo and Ni \cite{carrillo2009sharp}, Yokota \cite{yokota2012addendum}. It remains interesting whether this equality is true under more general assumptions. As an example, one may easily compute that on a compact Ricci flat static ancient flow, it holds that $\bar{W}(x_0)=\log{\bar{V}(x_0)}=-\infty$.
\\

In the special case of dimension three, because of the simplicity of the geometry, we are able to drop one of the additional assumptions.
\\

\begin{Main_Theorem_2} \label{Main_Theorem_2}
Let $(M^3,g(\tau))_{\tau\in[0,\infty)}$ be an ancient solution to the backward Ricci flow such that $\displaystyle\sup_{x\in M}R(x,0)<\infty$. Then the conclusion of Theorem \ref{Main_Theorem} and Theorem \ref{Main_Corollary} hold for $(M^3,g(\tau))_{\tau\in[0,\infty)}$.
\end{Main_Theorem_2}

\bigskip

From the estimates that we carry out in this paper, we are able to establish the following Gaussian upper and lower bounds for the conjugate heat kernel coupled with the backward Ricci flow.
\\

\begin{Gaussian_bound} \label{Gaussian_bound}
There exists $c>0$ and $C<\infty$ depending only on $\kappa$ and the dimension $n$, such that on an ancient backward Ricci flow $(M^n,g(\tau))_{\tau\in[0,\infty)}$ with bounded nonnegative curvature operator and that is $\kappa$-noncollapsed, the conjugate heat kernel $u(p,\tau)$ based at $(x_0,0)$ satisfies
\begin{eqnarray} \label{eq:gaussian}
\frac{c}{\tau^{\frac{n}{2}}}\exp\left(-\frac{C}{\tau}dist_\tau(p,p(\tau))^2\right)\leq u(p,\tau)\leq \frac{C}{\tau^{\frac{n}{2}}}\exp\left(-\frac{c}{\tau}dist_\tau(p,p(\tau))^2\right),
\end{eqnarray}
for all $\tau>0$ and $p(\tau)$ is chosen such that $l(p(\tau),\tau)\leq\frac{n}{2}$, where $l$ is the reduced distance based also at $(x_0,0)$; see (\ref{eq:l_def}) for the definition of $l$.
\end{Gaussian_bound}

\bigskip

The above Gaussian bounds (\ref{eq:gaussian}) are somewhat surprising. They asserts that the conjugate heat kernel stays like a Gaussian kernel, but its ``center'', instead of being fixed at the base point, is moving along with the ``center'' of the reduced distance. One may imagine a conjugate heat kernel on a Bryant soliton based at the origin, as $\tau$ evolves to infinity, $p(\tau)$ will, of necessity, be drifted away to infinity, since along those $p(\tau)$'s one must eventually get a cylinder, which is the only possible asymptotic Ricci shrinker for the Bryant soliton. We hope this phenomenon will provide a better understanding of the relationship between two fundamental tools in the Ricci flow, the entropy and the reduced geometry.
\\

This paper is organized as follows. In section 2 we give all the definitions and collect some well-known results upon which our arguments are built, most of those results are due to Perelman. In section 3 we establish the estimates for the conjugate heat kernel, which turns out to share many properties with the reduced distance. Corollary \ref{Gaussian_bound} is proved at the end of section 3. In section 4 we use the estimates to prove Theorem \ref{Main_Theorem}, Theorem \ref{Main_Corollary}, and Theorem \ref{Main_Theorem_2}.

\section{Preliminaries}
\newtheorem{blow_down_shrinker}{Proposition}[section]
\newtheorem{basic_estimates}[blow_down_shrinker]{Lemma}
\newtheorem{basic_estimates_2}[blow_down_shrinker]{Lemma}
\newtheorem{distance_distortion}[blow_down_shrinker]{Lemma}
\newtheorem{Gaussian_concentration}[blow_down_shrinker]{Proposition}

In this section we collect some known results that are useful for our arguments, most of which are due to Perelman \cite{perelman2002entropy} and \cite{perelman2003ricci}, one can find detailed treatment in Morgan and Tian \cite{morgan2007ricci}. Unless otherwise specified, we always consider $(M,g(\tau))_{\tau\in[0,\infty)}$, an ancient $\kappa$-noncollapsed solution to the backward Ricci flow with bounded nonnegative curvature operator. Let $(x_0,0)$ be a fixed point in space-time.
\\

Let
\begin{eqnarray} \label{eq:def_kernel}
u(p,\tau)=\frac{1}{(4\pi\tau)^{\frac{n}{2}}}e^{-f(p,\tau)}
\end{eqnarray}
be the fundamental solution to the conjugate heat equation coupled with the backward Ricci flow
\begin{eqnarray*}
\frac{\partial}{\partial\tau}g&=&2Ric,
\\
\frac{\partial}{\partial\tau}u&=&\Delta u-Ru,
\end{eqnarray*}
where $(x_0,0)$ is the base point such that $\displaystyle\lim_{\tau\rightarrow 0+}u(p,\tau)=\delta_{x_0}(p)$. Then Perelman's \textit{pointed entropy} is defined as
\begin{eqnarray}
W_{x_0}(\tau)&=&\int_M\left\{\tau(|\nabla f|^2+R)+f-n\right\}\frac{1}{(4\pi\tau)^{\frac{n}{2}}}e^{-f}dg(\tau) \label{eq:def_1}
\\
&=&\int_M\left\{\tau(2\Delta f-|\nabla f|^2+R)+f-n\right\}\frac{1}{(4\pi\tau)^{\frac{n}{2}}}e^{-f}dg(\tau). \label{eq:def_2}
\end{eqnarray}
It is known that
\begin{eqnarray*}
\lim_{\tau\rightarrow 0+}W_{x_0}(\tau)=0.
\end{eqnarray*}
In \cite{zhang2017entropy} we have defined the \textit{asymptotic entropy} as
\begin{eqnarray} \label{eq:def_W}
\bar{W}(x_0)=\lim_{\tau\rightarrow\infty}W_{x_0}(\tau),
\end{eqnarray}
where the limit always exists (and could possibly be $-\infty$) since $W_{x_0}(\tau)$ is monotonically decreasing in $\tau$
\begin{eqnarray} \label{eq:monotone}
\frac{d}{d\tau}W_{x_0}(\tau)=-\int_M2\tau\left|Ric+Hessf-\frac{1}{2\tau}g\right|^2\frac{1}{(4\pi\tau)^{\frac{n}{2}}}e^{-f}dg(\tau).
\end{eqnarray}
\\

Recall that Perelman's \textit{reduced distance} $l(p,\tau)$ based at $(x_0,0)$ is defined as
\begin{eqnarray} \label{eq:l_def}
l(p,\tau)=\frac{1}{2\sqrt{\tau}}\inf_{\gamma}\left(\int_0^{\tau}\sqrt{s}\left(R(\gamma(s),s)+|\dot{\gamma}(s)|^2_{g(s)}\right)ds\right),
\end{eqnarray}
where $\tau>0$ and the infimum is taken among all the curves $\gamma(s):[0,\tau]\rightarrow M$ such that $\gamma(0)=x_0$ and $\gamma(\tau)=p$. The minimizing curve is called the shortest $\mathcal{L}$-geodesic. Perelman's \textit{reduced volume} is defined as
\begin{eqnarray}
V_{x_0}(\tau)=\int_M(4\pi\tau)^{-\frac{n}{2}}e^{-l(p,\tau)}dg(p,\tau).
\end{eqnarray}
Taking the limit of the reduced volume as $\tau\rightarrow\infty$, we obtain the \textit{asymptotic reduce volume}

\begin{eqnarray}\label{def:asymp_rvol}
\bar{V}(x_0)=\lim_{\tau\rightarrow\infty}V_{x_0}(\tau),
\end{eqnarray}
where the limit always exists since $V_{x_0}(\tau)$ is monotonically decreasing in $\tau$.
\\

It is known from Perelman \cite{perelman2002entropy} that $\inf_{p\in M}l(p,\tau)\leq \frac{n}{2}$ for any $\tau>0$. Let $\{\tau_k\}_{k=1}^\infty$ and $\{p_k\}_{k=1}^\infty$ be such that $\tau_k\nearrow\infty$ and $l(p_k,2\tau_k)\leq\frac{n}{2}$. Perelman proved the existence of the asymptotic Ricci shrinker, as a limit of the blow-down sequence centered at $\{p_k\}_{k=1}^\infty$ and with respect to the scaling factors $\{\tau_k^{-1}\}_{k=1}^\infty$.
\\

\begin{blow_down_shrinker} [Perelman \cite{perelman2002entropy}] \label{blow_down_shrinker}
By passing to a (not relabelled) subsequence, $\{(M^n,\\\tau_k^{-1}g(\tau\tau_k),(p_k,1))_{\tau\in[\frac{1}{2},2]}\}_{k=1}^\infty$ converges in the pointed smooth Cheeger-Gromov sense to the canonical form of a shrinking gradient Ricci soliton $(M_\infty,g_\infty(\tau),l_\infty(\tau))_{\tau\in[\frac{1}{2},2]}$, which is called the asymptotic Ricci shrinker. The asymptotic Ricci shrinker is normalized in the way that
\begin{eqnarray} \label{eq:normalize_0}
2\Delta l_\infty-|\nabla l_\infty|^2+R+\frac{l_\infty-n}{\tau}\equiv 0,
\end{eqnarray}
where $l_\infty(\tau)$ is the limit of $l_k(\tau)=l(\tau\tau_k)$.
\end{blow_down_shrinker}

One may refer to Theorem 9.11 in \cite{morgan2007ricci} for a detailed proof of the existence of the asymptotic Ricci shrinker, and Proposition 9.20 in the same book for the validity of (\ref{eq:normalize_0}). Notice that our choice of the space-time base points $(p_k,\tau_k)$ is somewhat different from theirs, but their proof still works. It is not hard to see that the essential condition needed for the choice of the space-time base points is that $l(p_k,\tau_k)$ is bounded uniformly in $k$.
\\

\begin{basic_estimates} [Perelman \cite{perelman2002entropy}]\label{basic_estimates}
For $f$ and $l$ defined in (\ref{eq:def_kernel}) and (\ref{eq:l_def}), respectively, and both based at $(x_0,0)$, the following holds.
\begin{enumerate}[(a)]
\item $f(p,\tau)\leq l(p,\tau)$, for all $(p,\tau)\in M\times(0,\infty)$.

\item There exists a constant $C<\infty$ such that
  \begin{eqnarray}
  |\nabla l|^2+R &\leq& \frac{Cl}{\tau}, \label{eq:rdistance_1}
  \\
  \frac{\partial l}{\partial \tau}+\frac{Cl}{\tau}&\geq& R. \label{eq:r_distance_2}
  \end{eqnarray}
  In particular, we have that $l(p_k,\tau)\leq C$, for all $\tau\in[\frac{1}{2}\tau_k,2\tau_k]$ and for all $k$.

\item The reduce distance has quadratic growth. For all $p$, $q\in M$ and $\tau>0$, it holds that
  \begin{eqnarray} \label{eq:l_quadratic}
  -2l(q,\tau)-C+\frac{c}{\tau}dist_{\tau}(p,q)^2\leq l(p,\tau)\leq 2l(q,\tau)+\frac{C}{\tau}dist_{\tau}(p,q)^2,
  \end{eqnarray}
  where $c$ and $C$ are constants (see Lemma 9.25 in \cite{morgan2007ricci} for a detailed proof).
\end{enumerate}
\end{basic_estimates}

\bigskip

\begin{basic_estimates_2}\label{basic_estimates_2}

\begin{enumerate}[(a)]
  \item (Bounded curvature at bounded distance; see 11.7 in \cite{perelman2002entropy}) There exists a positive function $C(r):(0,\infty)\rightarrow(0,\infty)$ depending only on $\kappa$ and the dimension $n$, such that
  \begin{eqnarray*}
      R(q,\tau)\leq C(A)R(p,\tau),
  \end{eqnarray*}
  whenever $dist_\tau(p,q)\leq A R(p,\tau)^{-\frac{1}{2}}$ and for all $\tau\geq 0$.
  \item (Perelman's derivative estimates; see 1.5 in \cite{perelman2003ricci})There exists $C<\infty$ depending only on $\kappa$ and the dimension $n$, such that
  \begin{eqnarray}
  \left|\frac{\partial R}{\partial\tau}\right|&\leq& CR^2, \label{eq:deriv_1}
  \\
  |\nabla R|&\leq&CR^{\frac{3}{2}}, \label{eq:deriv_2}
  \end{eqnarray}
  hold on $M\times [0,\infty)$. More generally, we have
  \begin{eqnarray} \label{eq:deriv_3}
  \left|\frac{\partial^l}{\partial \tau^l}\nabla^m Rm\right|\leq\eta R^{1+l+\frac{m}{2}},
  \end{eqnarray}
  for any positive integers $l$ and $m$, where $\eta$ depends on $\kappa$, $n$, $l$, and $m$.
\end{enumerate}
\end{basic_estimates_2}

\bigskip

\begin{distance_distortion} [Perelman \cite{perelman2002entropy}]\label{distance_distortion}
On any backward Ricci flow, suppose $Ric(x,\tau_0)\leq (n-1)K$ for all $x\in B_{\tau_0}(x_0,r_0)\cup B_{\tau_0}(x_1,r_0)$. Then
\begin{eqnarray}
\frac{d}{d\tau}dist_\tau(x_0,x_1)\leq 2(n-1)(\frac{2}{3}Kr_0+r_0^{-1})\ \ \text{at}\ \tau=\tau_0.
\end{eqnarray}
\end{distance_distortion}

\bigskip

We conclude this section with the following Gaussian concentration inequality by Hein and Naber \cite{hein2014new}, which is very important to the pointwise estimate of the conjugate heat kernel, and is the starting point of our work. Notice that in \cite{hein2014new} they only require the Ricci flow to have bounded geometry on every time slice, which is guaranteed in our case by the bounded curvature and noncollapsing assumption. Their result is in the spirit of Davies \cite{davies1990heat}.
\\

\begin{Gaussian_concentration} \label{Gaussian_concentration}
Let $u(p,\tau)$ be the conjugate heat kernel based at $(x_0,0)$ as defined in (\ref{eq:def_kernel}) and $d\nu_\tau=u(\tau)dg(\tau)$ be a probability measure. Let $A$, $B$ be two subsets of $M$, then
\begin{eqnarray} \label{eq:gaussian_concentration}
\nu_\tau(A)\nu_\tau(B)\leq\exp\left(-\frac{1}{8\tau}dist_\tau(A,B)^2\right).
\end{eqnarray}
\end{Gaussian_concentration}

\bigskip

\section{Estimates on the conjugate heat kernel}
\newtheorem{local_lower_bound}{Lemma}[section]
\newtheorem{heatkernel_pointwise}[local_lower_bound]{Proposition}
\newtheorem{locally_uniform}[local_lower_bound]{Proposition}

We will first of all fix the notations for the rest of this paper. Unless otherwise indicated, we let $(M,g(\tau))_{\tau\in[0,\infty)}$ be an ancient $\kappa$-noncollapsed solution to the backward Ricci flow with bounded nonnegative curvature operator, on which we fix a point $(x_0,0)$. Let $u(p,\tau)$ and $l(p,\tau)$ be the fundamental solution to the conjugate heat equation and the reduced distance as defined in (\ref{eq:def_kernel}) and (\ref{eq:l_def}), respectively, both based at $(x_0,0)$. Let $\{p_k\}_{k=1}^\infty$ and $\{\tau_k\}_{k=1}^\infty$ be the sequences of points and times as described at the beginning of the previous section, that is, $\tau_k\nearrow\infty$ and $l(p_k,2\tau_k)\leq\frac{n}{2}$. We denote
\begin{eqnarray*}
g_k(\tau)&=&\tau_k^{-1}g(\tau\tau_k),
\\
R_k&=&R(g_k),
\\
l_k(p,\tau)&=&l(p,\tau\tau_k),
\\
u_k(p,\tau)&=&\tau_k^{\frac{n}{2}}u(p,\tau\tau_k):=\frac{1}{(4\pi\tau)^{\frac{n}{2}}}e^{-f_k(p,\tau)}.
\end{eqnarray*}
By the parabolic scaling, $u_k$ is a fundamental solution to the conjugate heat equation coupled with the backward Ricci flow $g_k(\tau)$. For most of the time we will use $C$ and $c$ to denote positive estimation constants, which could be different from line to line. We use the capital letter $C$ to denote those constants that are intuitively large, and the lower case letter $c$ to denote those constants that are intuitively small. We start with observing the following lower bound of $u_k$.
\\

\begin{local_lower_bound}[The lower bound] \label{local_lower_bound}
There exists a constant $C<\infty$ such that
\begin{eqnarray}\label{eq:lower_bound}
f_k(p,\tau)\leq C+Cdist_{g_k(\tau)}(p_k,p)^2,
\end{eqnarray}
for all $p\in M$ and $\tau\in[\frac{1}{2},2]$. In particular, $u_k(p,\tau)>c$ for all $p\in B_{g_{k}(\tau)}(p_k,1)$ and $\tau\in[\frac{1}{2},2]$, where $c>0$ is a constant.
\end{local_lower_bound}

The above lemma is an immediate consequence of Lemma \ref{basic_estimates}(a)(c), both of which are invariant under parabolic scaling.
\\

Next we apply Hein-Naber's Gaussian concentration inequality to obtain an integral estimate, and furthermore pointwise estimate of $u_k$.
\\

\begin{heatkernel_pointwise} \label{heatkernel_pointwise}
There exists $c>0$ and $C<\infty$ depending only on $\kappa$ and the dimension $n$ such that
\begin{eqnarray} \label{eq:HK_pointwise}
f_k(p,\tau)\geq-C-C\log\left(\max\{R_k(p,\tau),1\}\right)+cdist_{g_k(\tau)}(p,p_k)^2,
\end{eqnarray}
for all $p\in M$ and $\tau\in[\frac{3}{4},2]$.
\end{heatkernel_pointwise}

\begin{proof}
First of all, we use Proposition \ref{Gaussian_concentration} to obtain an integral estimate for $u_k$, that is, there exists $C<\infty$ depending only on $\kappa$ and the dimension $n$ such that for all $(p,\tau)$ with $\tau\in[\frac{1}{2},2]$ and for all $r\in(0,1]$, it holds that
\begin{eqnarray} \label{eq:oneslice_integral}
\int_{B_{g_k(\tau)}(p,r)}u_k(x,\tau)dg_k(x,\tau)\leq C\exp\left(-\frac{1}{32}dist_{g_k(\tau)}(p_k,p)^2\right).
\end{eqnarray}
(\ref{eq:oneslice_integral}) follows from applying (\ref{eq:gaussian_concentration}) to $A=B_{g_k(\tau)}(p_k,1)$ and $B=B_{g_k(\tau)}(p,r)$. Notice that by Lemma \ref{local_lower_bound} we have $u_k(x,\tau)\geq c$ for all $x\in B_{g_k(\tau)}(p_k,1)$, and that by (\ref{eq:rdistance_1}) and (\ref{eq:l_quadratic}) we have $R_k(x,\tau)\leq C$ for $x\in B_{g_k(\tau)}(p_k,1)$, hence $\operatorname{Vol}(B_{g_k(\tau)}(p_k,1))>c$ by the noncollapsing assumption. Therefore, $\nu(A)=\int_{B_{g_k(\tau)}(p_k,1)}u_k(x,\tau)dg_\tau(x)>c$ and (\ref{eq:oneslice_integral}) follows immediately. Notice that both (\ref{eq:rdistance_1}) and (\ref{eq:l_quadratic}) are invariant under parabolic scaling.
\\

Next, we will extend inequality (\ref{eq:oneslice_integral}) to a space-time cube, so as to apply the parabolic mean value inequality. We fix $\tau\in[\frac{3}{4},2]$ and $p\in M$, and define $r=\min\{1,R_k(p,\tau)^{-\frac{1}{2}}\}$. Applying Lemma \ref{basic_estimates_2} we can find $c_1\in(0,\frac{1}{4}]$ and $C_1<\infty$ such that $R_k(x,s)\leq C_1r^{-2}$, for all $x\in B_{g_k(s)}(p,r)$ and $s\in[\tau-c_1r^2,\tau]$. Moreover, applying Lemma \ref{distance_distortion} to the points $p$ and $p_k$, we have that
\begin{eqnarray} \label{eq:dist_estimate}
dist_{g_k(s)}(p_k,p)\geq dist_{g_k(\tau)}(p_k,p)-10nC_1\sqrt{c_1},
\end{eqnarray}
for all $s\in[\tau-c_1r^2,\tau]$. Integrating (\ref{eq:oneslice_integral}) in time and taking (\ref{eq:dist_estimate}) into account, we obtain the following integral estimate
\begin{eqnarray*}
\int_{\tau-c_1r^2}^\tau\int_{B_{g_k(s)}(p,r)}u_k(x,s)dg_k(x,s)ds\leq c_1Cr^2\exp\left(-\frac{1}{64}dist_{g_k(\tau)}(p_k,p)^2\right).
\end{eqnarray*}
By a local distortion argument, we can always choose $c_2>0$ small enough and $C_2<\infty$ large enough, depending on all the previous estimation constants, such that
\begin{eqnarray*}
B_{g_k(\tau)}(p,c_2r)\subset B_{g_k(s)}(p,r)&,& \text{for all}\ s\in[\tau-(c_2r)^2,\tau],
\\
R_k(x,s)\leq C_2r^{-2}&,& \text{for all}\ (x,s)\in B_{g_k(\tau)}(p,c_2r)\times [\tau-(c_2r)^2,\tau],
\\
C_2^{-1}g_k(s)\leq g_k(\tau) \leq C_2 g_k(s)&,& \text{in}\ B_{g_k(\tau)}(p,c_2r)\times [\tau-(c_2r)^2,\tau],
\end{eqnarray*}
and that
\begin{eqnarray*}
\int_{\tau-(c_2r)^2}^\tau\int_{B_{g_k(\tau)}(p,c_2r)}u_k(x,s)dg_k(x,\tau)ds\leq C_2r^2\exp\left(-\frac{1}{64}dist_{g_k(\tau)}(p_k,p)^2\right).
\end{eqnarray*}
Now we are in the position to apply the parabolic mean value inequality (c.f. Theorem 25.12 of \cite{chow2010ricci}). To keep track with the curvature condition, we need only to do parabolic scaling in our scenario. Let $\tilde{u}(s)=u_k(\tau+sr^2)$ and $\tilde{g}(s)=r^{-2}g_k(\tau+sr^2)$, then we have
\begin{eqnarray*}
&&\int_{-(c_2)^2}^0\int_{B_{\tilde{g}(0)}(p,c_2)}\tilde{u}(x,s)d\tilde{g}(x,0)ds
\\
&&= \frac{1}{r^{n+2}}\int_{\tau-(c_2r)^2}^\tau\int_{B_{g_k(\tau)}(p,c_2r)}u_k(x,s)dg_k(x,\tau)ds \leq \frac{C_2}{r^n}\exp\left(-\frac{1}{64}dist_{g_k(\tau)}(p_k,p)^2\right),
\end{eqnarray*}
with the curvature bound
\begin{eqnarray*}
\tilde{R}(x,s)\leq C_2,\ \text{for all}\ (x,s)\in B_{\tilde{g}(0)}(p,c_2)\times [-(c_2)^2,0].
\end{eqnarray*}
Therefore we have
\begin{eqnarray*}
u_k(p,\tau)&=&\tilde{u}(p,0)\leq \frac{C}{(\frac{c_2}{4})^2\operatorname{Vol}(B_{\tilde{g}(0)}(p,\frac{c_2}{4}))}\cdot \frac{C_2}{r^n}\exp\left(-\frac{1}{64}dist_{g_k(\tau)}(p_k,p)^2\right) \\
&\leq& \frac{C}{r^n}\exp\left(-\frac{1}{64}dist_{g_k(\tau)}(p_k,p)^2\right),
\end{eqnarray*}
where in the last inequality we have again used the volume lower bound implied by the noncollapsing condition. Taking into account the definition of $r=\min\{1,R_k(p,\tau)^{-\frac{1}{2}}\}$ we obtain the conclusion.
\end{proof}

\bigskip

Once we have the above pointwise estimate for the conjugate heat kernel, we can apply the curvature growth condition given by (\ref{eq:rdistance_1}) to obtain the following nice quadratic growth estimate for $f_k$'s, and furthermore obtain a gradient growth estimate by using a localized gradient estimate of the conjugate heat equation.
\\

\begin{locally_uniform}\label{locally_uniform}
The following holds for every element of the sequence of backward Ricci flows $\{(M,g_k(\tau),(p_k,1))_{\tau\in[\frac{1}{2},2]}\}_{k=1}^\infty$.
\begin{enumerate}[(a)]
  \item There exists $C<\infty$ such that
  \begin{eqnarray} \label{eq:L_U_estimate_1}
  R_k(p,\tau)\leq C+Cdist_{g_k(\tau)}(p_k,p)^2,
  \end{eqnarray}
  for all $p\in M$ and $\tau\in[\frac{1}{2},2]$.
  \item There exists $C<\infty$ and $c>0$ depending only on $\kappa$ and the dimension $n$, such that
  \begin{eqnarray} \label{eq:L_U_estimate_2}
  f_k(p,\tau)\geq -C +cdist_{g_k(\tau)}(p_k,p)^2,
  \end{eqnarray}
  for all $p\in M$ and $\tau\in[\frac{3}{4},2]$.
  \item There exists $C<\infty$ depending only on $\kappa$ and the dimension $n$, such that
  \begin{eqnarray}\label{eq:L_U_estimate_3}
  |\nabla f_k(p,\tau)|^2\leq C+Cdist_{g_k(\tau)}(p_k,p)^7,
  \end{eqnarray}
  for all $p\in M$ and $\tau\in[1,2]$.
\end{enumerate}
In particular, every constant in these estimates do not depend on $k$.
\end{locally_uniform}

\begin{proof}
\begin{enumerate}[(a)]
  \item Scaling (\ref{eq:rdistance_1}) and the second inequality of (\ref{eq:l_quadratic}) by the factor $\tau_k^{-1}$, we have
  \begin{eqnarray*}
  R_k\leq|\nabla l_k|^2+R_k\leq\frac{Cl_k}{\tau}
  \end{eqnarray*}
  and
  \begin{eqnarray*}
  l_k(p,\tau)\leq C+\frac{C}{\tau}dist_{g_k(\tau)}(p_k,p)^2,
  \end{eqnarray*}
  where we have let $q=p_k$ in (\ref{eq:l_quadratic}) and used Lemma \ref{basic_estimates}(b) that $l_k(p_k,\tau)\leq C$ for all $\tau\in[\frac{1}{2},2]$. The conclusion of (a) follows immediately.

  \item Inserting (\ref{eq:L_U_estimate_1}) into (\ref{eq:HK_pointwise}) we have that
     \begin{eqnarray*}
     f_k(p,\tau)&\geq&-C-C\log\left(C+Cdist_{g_k(\tau)}(p_k,p)^2\right)+cdist_{g_k(\tau)}(p,p_k)^2
     \\
     &\geq&-C+\frac{c}{2}dist_{g_k(\tau)}(p,p_k)^2
     \\
     &&+\left(\frac{c}{2}dist_{g_k(\tau)}(p,p_k)^2-C\log\left(C+Cdist_{g_k(\tau)}(p_k,p)^2\right)\right).
     \end{eqnarray*}
     The conclusion of (b) follows from the simple fact that the function $\displaystyle \frac{c}{2}x^2-C\log(C+Cx^2)$ is bounded from below.

  \item We apply Theorem 10 of \cite{ecker2008local}, where they have shown if on $\Omega(2A)=\bigcup_{\tau\in[0,\bar{t}]}B_{g(\tau)}(q,2A)$ the following bounds hold
  \begin{eqnarray*}
  |Ric|\leq K_1\ \text{and}\ |\nabla R|\leq K_2,
  \end{eqnarray*}
  and $u$ is a positive solution to the conjugate heat equation with $u\leq J$, then it holds that
  \begin{eqnarray*}
  \frac{|\nabla u|^2}{u^2}\leq\left(1+\log\frac{J}{u}\right)^2\left(\frac{1}{\tau}+C_1K_1+\sqrt{K_2}+K_2+\frac{C_1\sqrt{K_2}A\coth(\sqrt{K_2}A)+C_2}{A^2}\right)
  \end{eqnarray*}
  on $\Omega(A)$ and for $\tau\in(0,\bar{t}]$, where $C_1$ and $C_2$ are constants. In order to apply this theorem to each $u_k$ on $[\frac{3}{4},2]$, we need to check every bound that is needed. By (\ref{eq:L_U_estimate_2}) we have a uniform upper bound $J$ for $u_k$. By (\ref{eq:L_U_estimate_1}) and (\ref{eq:deriv_2}), we may take $K_1=C+CA^2$ and $K_2=C+CA^3$. Moreover, (\ref{eq:lower_bound}) provides an upper bound for $\log\frac{J}{u}\leq C +CA^2$. Hence we have
  \begin{eqnarray*}
  \frac{|\nabla u_k|^2}{u_k^2}(p,\tau)\leq C(1+A^2)^2\left(1+A^\frac{3}{2}+A^2+A^3+\frac{A(1+A^{\frac{3}{2}})\coth(CA(1+A^{\frac{3}{2}}))+1}{A^2}\right)
  \end{eqnarray*}
  for all $\displaystyle p\in B_{g_{k}(\tau)}(p_k,A)$ and $\tau\in[1,2]$, whence follows the result.
\end{enumerate}
\end{proof}

\bigskip

\begin{proof}[Proof of Corollary \ref{Gaussian_bound}]
Notice that in all the estimates above, we may indeed take $\tau_k$ to be any positive number. Taking $\tau=2$ in the estimates (\ref{eq:lower_bound}) and (\ref{eq:L_U_estimate_2}), and scaling them back by the factor $\tau_k$, we obtain the conclusion.
\end{proof}

\section{Proof of the main theorem}

\newtheorem{No_loss_of_entropy}{Proposition}[section]

In Proposition \ref{locally_uniform} we have derived a uniform upper bound for $u_k$'s on the interval $[1,2]$. Combining (\ref{eq:deriv_3}) and (\ref{eq:L_U_estimate_1}) we obtain locally uniform estimates for all the derivatives of the curvatures on $[1,2]$. By standard parabolic regularity theory, we have locally uniform estimates for the solutions and all of their derivatives. Hence by passing to a subsequence, $\{u_k\}_{k=1}^\infty$ converges locally smoothly to a solution $\displaystyle u_\infty:=\frac{1}{(4\pi\tau)^{-\frac{n}{2}}}e^{-f_\infty}$ to the conjugate heat equation on the asymptotic Ricci shrinker, with $f_k\rightarrow f_\infty$ locally smoothly. Here we would like to remark that the parabolic regularity theory actually ensures the smooth convergence of $\{u_k\}_{k=1}^\infty$ on any space-time compact subset of $M\times[1,2)$, however, it is not hard to extend all those estimates we obtained in the previous section to an interval larger than $[1,2]$, say $[\frac{7}{8},4]$. Therefore, to keep our notations concise, we simply assume that the smooth convergence of $\{u_k\}_{k=1}^\infty$ happens on any compact subset of $M\times[1,2]$. The curvature growth condition (\ref{eq:L_U_estimate_1}), the quadratic upper and lower bounds (\ref{eq:lower_bound}) and (\ref{eq:L_U_estimate_2}) of $f_k$, and the derivative growth bound (\ref{eq:L_U_estimate_3}) carries to $R_\infty$ and $f_\infty$ on $M_\infty\times[1,2]$:
\begin{eqnarray}
R_\infty(p,\tau) &\leq& C+Cdist_{g_\infty(\tau)}(p_\infty,p)^2, \label{eq:limit_curvature}
\\
  -C+cdist_{g_\infty(\tau)}(p_\infty,p)^2&\leq& f_\infty(p,\tau)\leq C+Cdist_{g_\infty(\tau)}(p_\infty,p)^2, \label{eq:limit_quadratic}
\\
  |\nabla f_\infty(p,\tau)|^2&\leq& C+Cdist_{g_\infty(\tau)}(p_\infty,p)^7, \label{eq:deriv_growth}
\end{eqnarray}
for all $p\in M_\infty$ and $\tau\in[1,2]$. In particular, (\ref{eq:limit_quadratic}) implies $u_\infty$ is not a zero solution.
\\

We denote
\begin{eqnarray*}
W_k(\tau)&=&\int_{M}\left\{\tau(|\nabla f_k|^2+R_k)+f_k-n\right\}\frac{1}{(4\pi\tau)^{\frac{n}{2}}}e^{-f_k}dg_k(\tau),
\\
W_\infty(\tau)&=&\int_{M_\infty}\left\{\tau(|\nabla f_\infty|^2+R_\infty)+f_\infty-n\right\}\frac{1}{(4\pi\tau)^{\frac{n}{2}}}e^{-f_\infty}dg_\infty(\tau).
\end{eqnarray*}

\bigskip

Theorem \ref{Main_Theorem} is implied by the following no loss of entropy proposition as well as Theorem 1.2 in \cite{zhang2017entropy}.
\\

\begin{No_loss_of_entropy}
$W_\infty(\tau)=W_\infty>-\infty$ is a constant, and $\displaystyle \lim_{\tau\rightarrow\infty}W_{x_0}(\tau)=W_\infty$. Moreover, $f_\infty$ is a potential function of the asymptotic Ricci shrinker.
\end{No_loss_of_entropy}

\begin{proof}

\newtheorem{Claim_1}{Claim}
\newtheorem{Claim_2}[Claim_1]{Claim}
\newtheorem{Claim_3}[Claim_1]{Claim}

We split the proof into several steps, first we show that $\displaystyle \liminf_{k\rightarrow\infty}W_k(\tau)\geq W_\infty(\tau)>-\infty$, and then we show that $\displaystyle \limsup_{k\rightarrow\infty}W_k(\tau)\leq W_\infty(\tau)$, finally we show that $W_\infty(\tau)$ is a constant.

\begin{Claim_1}
There exists a $A_0<\infty$ independent of $k$, such that
\begin{eqnarray}
\{\tau(|\nabla f_k|^2+R_k)+f_k-n\}(p,\tau)&\geq&0 \label{eq:claim_1_1}
\\
\{\tau(|\nabla f_\infty|^2+R_\infty)+f_\infty-n\}(p,\tau)&\geq& 0 \label{eq:claim_1_2}
\end{eqnarray}
whenever $dist_{g_k(\tau)}(p_k,p)\geq A_0$, $dist_{g_\infty(\tau)}(p_\infty,p)\geq A_0$, and for all $\tau\in[1,2]$.
\end{Claim_1}

Claim 1 follows immediately from the quadratic growth estimates (\ref{eq:L_U_estimate_2}) and (\ref{eq:limit_quadratic}) of $f_k$ and $f_\infty$, respectively.
\\

\begin{Claim_2}
$\displaystyle \liminf_{k\rightarrow\infty}W_k(\tau)\geq W_\infty(\tau)>-\infty$
\end{Claim_2}

\begin{proof}[Proof of the claim]
By the local convergence of $\{f_k(p,\tau)\}_{k=1}^\infty$ and $\{u_k(p,\tau)\}_{k=1}^\infty$, we have that for all $A>A_0$ and $\tau\in[1,2]$, the following holds
\begin{eqnarray*}
\liminf_{k\rightarrow\infty} W_k(\tau)&\geq& \lim_{k\rightarrow\infty}\int_{B_{g_k(\tau)}(p_k,A)}\left\{\tau(|\nabla f_k|^2+R_k)+f_k-n\right\}\frac{1}{(4\pi\tau)^{\frac{n}{2}}}e^{-f_k}dg_k(\tau)
\\
&=&\int_{B_{g_\infty(\tau)}(p_\infty,A)}\left\{\tau(|\nabla f_\infty|^2+R_\infty)+f_\infty-n\right\}\frac{1}{(4\pi\tau)^{\frac{n}{2}}}e^{-f_\infty}dg_\infty(\tau)
\\
&\geq&\int_{B_{g_\infty(\tau)}(p_\infty,A_0)}\left\{\tau(|\nabla f_\infty|^2+R_\infty)+f_\infty-n\right\}\frac{1}{(4\pi\tau)^{\frac{n}{2}}}e^{-f_\infty}dg_\infty(\tau)
\\
&>&-\infty.
\end{eqnarray*}
Taking $A\rightarrow\infty$ in the above inequality, we have $\displaystyle \liminf_{k\rightarrow\infty}W_k(\tau)\geq W_\infty(\tau)>-\infty$.
\end{proof}

\smallskip

\begin{Claim_3}
 $\displaystyle \limsup_{k\rightarrow\infty}W_k(\tau)\leq W_\infty(\tau)$.
\end{Claim_3}

\begin{proof}[Proof of the claim]
We use the equivalent definition (\ref{eq:def_2}) of Perelman's entropy. By Perelman's differential Harnack inequality (section 9 in \cite{perelman2002entropy}), we have that the integrand in (\ref{eq:def_2}) is nonpositive
\begin{eqnarray*}
\tau(2\Delta f_k-|\nabla f_k|^2+R_k)+f_k-n\leq 0,
\end{eqnarray*}
for all $\tau\in[1,2]$ and for all $k$, obviously this inequality carries to $f_\infty$
\begin{eqnarray*}
\tau(2\Delta f_\infty-|\nabla f_\infty|^2+R_\infty)+f_\infty-n\leq 0
\end{eqnarray*}
because of the locally smooth convergence.
\\

Let $A<\infty$ and $0\leq\phi_k\leq 1$ be cut-off functions such that $\phi_k=1$ on $B_{g_k(\tau)}(p_k,A)$, $\phi_k=0$ outside $B_{g_k(\tau)}(p_k,A+1)$, $|\nabla \phi_k|\leq 2$, and that $\phi_k\rightarrow\phi_\infty$ uniformly, where $\phi_\infty$ is a cut-off function on $M_\infty$ with the same properties.
\begin{eqnarray*}
&&\limsup_{k\rightarrow\infty}W_k(\tau)
\\
&&\leq\lim_{k\rightarrow\infty}\int_M\phi_k^2\left\{\tau(2\Delta f_k-|\nabla f_k|^2+R_k)+f_k-n\right\}\frac{1}{(4\pi\tau)^{\frac{n}{2}}}e^{-f_k}dg_k(\tau)
\\
&&=\int_{M_\infty}\phi_\infty^2\left\{\tau(2\Delta f_\infty-|\nabla f_\infty|^2+R_\infty)+f_\infty-n\right\}\frac{1}{(4\pi\tau)^{\frac{n}{2}}}e^{-f_\infty}dg_\infty(\tau).
\end{eqnarray*}
When we take $A\rightarrow\infty$ in the last formula, except for $\Delta f_\infty$ term we have sufficient bounds (\ref{eq:limit_curvature}), (\ref{eq:limit_quadratic}), and (\ref{eq:deriv_growth}). To handle $\Delta f_\infty$ term, we first notice that (\ref{eq:limit_quadratic}) and (\ref{eq:deriv_growth}) guarantees that the following integration by parts formula holds
\begin{eqnarray*}
\int_{M_\infty}\Delta f_\infty e^{-f_\infty}dg_\infty=\int_{M_\infty}|\nabla f_\infty|^2 e^{-f_\infty}dg_\infty.
\end{eqnarray*}
Furthermore we have
\begin{eqnarray*}
\int_{M_\infty}\phi_\infty^2\Delta f_\infty e^{-f_\infty}dg_\infty=\int_{M_\infty}\phi_\infty^2|\nabla f_\infty|^2 e^{-f_\infty}dg_\infty-2\int_{M_\infty}\phi_\infty\langle\nabla\phi_\infty\nabla f_\infty\rangle e^{-f_\infty}dg_\infty,
\end{eqnarray*}
therefore
\begin{eqnarray*}
&&\int_{M_\infty}(1-\phi_\infty^2)\Delta f_\infty e^{-f_\infty}dg_\infty
\\
&&=\int_{M_\infty}(1-\phi_\infty^2)|\nabla f_\infty|^2 e^{-f_\infty}dg_\infty + 2\int_{M_\infty}\phi_\infty\langle\nabla\phi_\infty\nabla f_\infty\rangle e^{-f_\infty}dg_\infty.
\end{eqnarray*}
Taking $A\rightarrow\infty$ we have $\int_{M_\infty}\phi_\infty^2\Delta f_\infty e^{-f_\infty}dg_\infty\rightarrow\int_{M_\infty}\Delta f_\infty e^{-f_\infty}dg_\infty$, hence we have prove the claim.
\end{proof}

Finally, since we have established $\displaystyle \lim_{k\rightarrow\infty}W_k(\tau)= W_\infty(\tau)$ for all $\tau\in[1,2]$. It follows from the monotonicity of $W_{x_0}(\tau)$ that for all $\tau\in[1,2]$, $\displaystyle\lim_{k\rightarrow\infty}W_k(\tau)=\lim_{s\rightarrow\infty}W_{x_0}(s)$, and hence $\displaystyle W_\infty(\tau)\equiv\lim_{s\rightarrow\infty}W_{x_0}(s)$ is a constant. Finally, the estimates (\ref{eq:limit_quadratic}) and (\ref{eq:deriv_growth}) justify the integration by parts in the derivation of entropy monotonicity formula (\ref{eq:monotone}) for $W_\infty(\tau)$, it follows that $f_\infty$ is indeed a potential function of the asymptotic Ricci shrinker $(M_\infty,g_\infty)$.

\end{proof}

\textbf{Remark.} Here we would like to remind the reader of the special care we have taken in the proof of the above proposition. One needs to notice that Perelman's differential Harnack inequality is proved only for fundamental solutions to the conjugate heat equation under the assumption of bounded curvature, neither of which is the case for $f_\infty$, wherefore the differential Harnack inequality for $f_\infty$ that we applied is due to locally smooth convergence. Moreover, we do not know the sign or an exact growth estimate of the term $\Delta f_\infty$; that is why we used a careful cut-off argument to show that outside a large ball its contribution to the integral is negligible.

\bigskip

\begin{proof}[Proof of Theorem \ref{Main_Corollary}]
\newtheorem{Claim_2_1}{Claim}
In the case that $(M,g(\tau))_{\tau\in[0,\infty)}$ is collapsed, both sides of (\ref{eq:equivalence}) are $-\infty$, where $\bar{W}(x_0)=-\infty$ follows from Theorem \ref{Main_Theorem} and $\log\bar{V}(x_0)=-\infty$ follows from Perelman \cite{perelman2002entropy} and Yokota \cite{yokota2009perelman}.  Hence we only consider the noncollapsing case.
\\

First of all, we establish the fact that there is no integration loss when taking limits of $u_k$ and $(4\pi\tau)^{-\frac{n}{2}}e^{-l_k(p,\tau)}$.
\begin{Claim_2_1}
\begin{eqnarray}
\int_{M_\infty}u_\infty(p,\tau)dg_\infty(p,\tau)&\equiv& 1, \label{eq:normalize_1}
\\
\int_{M_\infty}\frac{1}{(4\pi\tau)^\frac{n}{2}}e^{-l_\infty(p,\tau)}dg_\infty(p,\tau)&\equiv&\bar{V}(x_0), \label{eq:normalize_2}
\end{eqnarray}
for all $\tau\in[1,2]$.
\end{Claim_2_1}

\begin{proof}[Proof of the claim]
To prove the first equality, we apply the quadratic growth estimate of $f_k$ (\ref{eq:L_U_estimate_2}) and the Bishop-Gromov volume comparison to see that the contribution of $u_k$ to the integral outside a large ball is negligible. To wit, for any $\varepsilon>0$, there exists $A_0<\infty$, such that for all $A>A_0$ and for all $k$, it holds that
\begin{eqnarray*}
1-\varepsilon<\int_{B_{g_k}(p_k,A)}u_k(\tau)dg_k(\tau)\leq 1.
\end{eqnarray*}
The conclusion follows from first taking $k\rightarrow\infty$, then taking $A\rightarrow\infty$, and finally taking $\varepsilon\rightarrow 0$. The second equality follows from the same argument.
\end{proof}

\smallskip

Now we study the Ricci shrinker potentials $f_\infty$ and $l_\infty$ with $\tau$ fixed at $1$. By Carrillo and Ni \cite{carrillo2009sharp}, we know that $\bar{W}(x_0)=W_\infty(1)=\mu(g_\infty,1)$, since $f_\infty$ is normalized in the way of (\ref{eq:normalize_1}). Moreover, if we define $\tilde{l}_\infty=l_\infty+\log\bar{V}(x_0)$, by (\ref{eq:normalize_2}) we have $\int_{M_\infty}(4\pi)^{-\frac{n}{2}}e^{-\tilde{l}_\infty}=1$ and $\tilde{l}_\infty$ is normalized in the same way as $f_\infty$. Therefore, by applying Carrillo and Ni again we have
\begin{eqnarray*}
\mu(g_\infty,1)&=&\int_{M_\infty}\{2\Delta \tilde{f}_\infty-|\nabla\tilde{l}_\infty|^2+R_\infty+\tilde{l}_\infty-n\}(4\pi)^{-\frac{n}{2}}e^{-\tilde{l}_\infty}dg_\infty
\\
&=&\log\bar{V}(x_0),
\end{eqnarray*}
where in the last equality we have used (\ref{eq:normalize_0}) and the definition of $\tilde{l}_\infty$. Hence we have proved the corollary.

\end{proof}

\bigskip

\begin{proof} [Proof of Theorem \ref{Main_Theorem_2}]
\newtheorem{Claim_3_1}{Claim}
\newtheorem{Claim_3_2}[Claim_3_1]{Claim}

In the case where the manifold is compact, we always have bounded geometry on one time slice. Let $(M^3,g(\tau))_{\tau\in[0,\infty)}$ be a noncompact ancient solution to the backward Ricci flow with bounded curvature, then it either has positive curvature everywhere or, by Hamilton's strong maximum principle \cite{hamilton1986four}, splits locally. In the first case, by a Theorem of Gromov and Meyer, there is a positive lower bound for the injectivity radius of $(M^3,g(0))$ (c.f. Theorem B.65 in \cite{chow2004ricci}), hence Theorem \ref{Main_Theorem} is applicable. In the second case, the universal cover of $(M^3,g(\tau))$ must be a line cross a two dimensional ancient solution with bounded curvature, which could only be the sphere, the cigar soliton, or the King-Rosenau solution; see Chu \cite{chu2007type}, Daskalopoulos and Sesum \cite{daskalopoulos2006eternal}, Daskalopoulos, Hamilton, and Sesum \cite{daskalopoulos2012classification}. We denote the cigar soliton as $\mathcal{C}^2$ and the King-Rosenau solution as $\mathcal{K}^2$. Since $\mathbb{S}^2\times\mathbb{R}$, $\mathcal{C}^2\times\mathbb{R}$, and $\mathcal{K}^2\times\mathbb{R}$ all have bounded geometry, we need to consider their quotients.
\\

\begin{Claim_3_1} \footnote{The author would like to thank Professor Bennett Chow for telling him about this argument.}
Let $(\mathcal{C}^2\times\mathbb{R})/\Gamma$ be any nontrivial isometric quotient of $\mathcal{C}^2\times\mathbb{R}$, then $\Gamma$ is isomorphic to $\mathbb{Z}$ and is generated by the composition of a translation in $\mathbb{R}$ and a rotation of $\mathcal{C}^2$.
\end{Claim_3_1}

\begin{proof}[Proof of the claim]
Denote the projection map to the $\mathbb{R}$ factor of $\mathcal{C}^2\times\mathbb{R}$ as $\pi$, the origin of the $\mathcal{C}^2$ factor and the $\mathbb{R}$ factor as $O$ and $0$, respectively. Let $\gamma\in\Gamma$ be any element that is not the identity. Since $\mathcal{C}^2$ is a surface of positive curvature, the only direction on $\mathcal{C}^2\times\mathbb{R}$ along which the Ricci curvature vanishes is the $\mathbb{R}$ direction, which must be preserved by the isometry $\gamma$. Therefore $d\gamma(\partial_r)=\pm\partial_r$, and it follows that either $d\gamma(\partial_r)\equiv \partial_r$ or $d\gamma(\partial_r)\equiv -\partial_r$ since $d\gamma$ is continuous. Furthermore we have that $\gamma(\mathcal{C}^2\times\{x\})$ is a slice of the direct product, for any $x\in\mathbb{R}$.
\\

We show that $\gamma$ does not flip the orientation of the $\mathbb{R}$ factor. Obviously we have $\gamma(\mathcal{C}^2\times\{x\})\neq \mathcal{C}^2\times\{x\}$, because otherwise $(O,x)$ is a fixed point of under the group action (notice the only isometries on $\mathcal{C}^2$ are the rotations with respect to $O$), contradiction. Suppose $d\gamma(\partial_r)\equiv -\partial_r$, we may assume without loss of generality that $x>0$ and $\gamma(\mathcal{C}^2\times\{x\})= \mathcal{C}^2\times\{-x\}$. Define $S=\{y\in\mathbb{R}:\pi\circ\gamma(p)=-y,\ \text{for all}\ p\in\mathcal{C}^2\times\{y\}\}$. It is easy to see that $S$ is closed and not empty. Since $d\gamma(\partial_r)\equiv-\partial_r$, $S$ is also open. Hence $S=\mathbb{R}$ and $0\in S$, which is a contradiction because $(O,0)$ is a fixed point under the group action. Hence $d\gamma(\partial_r)\equiv\partial_r$.
\\

Next we show that there exists $c\neq 0$, such that $\gamma(\mathcal{C}^2\times\{x\})=\mathcal{C}^2\times\{x+c\}$ for all $x\in\mathbb{R}$. We consider the function
\begin{eqnarray*}
f:\mathbb{R}&\rightarrow&\mathbb{R}
\\
x&\mapsto& \pi\circ\gamma(\mathcal{C}^2\times\{x\})-x.
\end{eqnarray*}
Since $d\gamma(\partial_r)\equiv \partial_r$, we have that $f$ is continuous and locally constant, and hence it is constant.
\\

Finally we prove that there is a generator of $\Gamma$. Since the action of $\Gamma$ is discrete, we can always find a $\gamma_0$, at which the minimum of $\displaystyle\min_{\gamma\neq id}|\pi\circ\gamma(\mathcal{C}^2\times\{0\})|$ is attained. Without loss of generality, assume $\pi\circ\gamma_0(\mathcal{C}^2\times\{0\})=c_0>0$, then for all $\gamma\in\Gamma$, there exists $n\in\mathbb{Z}$, such that $\pi\circ\gamma(\mathcal{C}^2\times\{0\})=nc_0$. Suppose this is not the case, then we can find $\gamma\in\Gamma$ and $n\in\mathbb{Z}$ such that $nc_0<\pi\circ\gamma(\mathcal{C}^2\times\{0\})<(n+1)c_0$, it follows that $0<\pi\circ(\gamma\cdot\gamma_0^{-n})(\mathcal{C}^2\times\{0\})<c_0$, a contradiction. Furthermore, if $\pi\circ\gamma(\mathcal{C}^2\times\{0\})=nc_0$, then we must have $\gamma=\gamma_0^n$, because otherwise we have $\gamma\cdot\gamma_0^{-n}\neq id$ but $\gamma\cdot\gamma_0^{-n}(\mathcal{C}^2\times\{0\})= \mathcal{C}^2\times\{0\}$, a contradiction.

\end{proof}

\bigskip

\begin{Claim_3_2}
Let $(\mathcal{K}^2\times \mathbb{R})/\Gamma$ be a quotient of $\mathcal{K}^2\times \mathbb{R}$, then either $\Gamma$ is finite or $(\mathcal{K}^2\times \mathbb{R})/\Gamma$ is closed; the same conclusion holds for $(\mathbb{S}^2\times \mathbb{R})/\Gamma$.
\end{Claim_3_2}

\begin{proof}[Proof of the claim]
We may assume $(\mathcal{K}^2\times \mathbb{R})/\Gamma$ is orientable, for otherwise we may consider its $\mathbb{Z}_2$ cover. By the same argument as in the proof of Claim 1, for all $\gamma\in\Gamma$, either $d\gamma(\partial_r)\equiv\partial_r$ or $d\gamma(\partial_r)\equiv-\partial_r$.
\\

If there exists $\gamma\neq id\in\Gamma$ such that $d\gamma(\partial_r)\equiv\partial_r$, then $\gamma(\mathcal{K}^2\times\{x\})\neq\mathcal{K}^2\times\{x\}$ for all $x\in\mathbb{R}$, since there is no orientation preserving isometry on $\mathcal{K}^2$ that has no fixed point, and $\gamma(\mathcal{K}^2\times\{x\})=\mathcal{K}^2\times\{x\}$ implies every point on $\mathcal{K}^2\times\{x\}$ is fixed. By the same argument as in the proof of Claim 1 we have that there exist $c\neq 0$ such that $\gamma(\mathcal{K}^2\times\{x\})=\mathcal{K}^2\times\{x+c\}$ for all $x\in\mathbb{R}$, it follows that $(\mathcal{K}^2\times \mathbb{R})/\Gamma$ is closed since $\mathcal{K}^2$ is closed.
\\

If for all $\gamma\neq id\in\Gamma$ it holds that $d\gamma(\partial_r)\equiv-\partial_r$, then $\Gamma=\{id,\gamma\}$ is finite. Suppose this is not the case, then we have $\gamma_1\neq\gamma_2\in\Gamma\setminus\{id\}$, $\gamma_1\cdot\gamma_2^{-1}\neq id$, and $d(\gamma_1\cdot\gamma_2)^{-1}(\partial_r)\equiv\partial_r$, a contradiction. The same proof carries to $(\mathbb{S}^2\times \mathbb{R})/\Gamma$ since there also exists no orientation preserving isometry on $\mathbb{S}^2$ without fixed point.

\end{proof}

\bigskip

If $M^3=(\mathcal{C}^2\times\mathbb{R})/\Gamma$, where $\Gamma$ is generated by $\gamma_0$ and $\gamma_0(\mathcal{C}^2\times\{x\})=\mathcal{C}^2\times\{x+c\}$ for all $x\in\mathbb{R}$ and some fixed $c\neq 0$. Then it is easy to see that $\operatorname{inj}(M^3)\geq\min\{\operatorname{inj}(\mathcal{C}^2),\frac{|c|}{2}\}>0$, hence there is a lower bound of the volume of the unit balls and Theorem \ref{Main_Theorem} is applicable.
\\

If $M^3=(\mathcal{K}^2\times\mathbb{R})/\Gamma$ or $(\mathbb{S}^2\times \mathbb{R})/\Gamma$, in the case $\Gamma$ is infinite, we have $M^3$ is closed and there is a lower bound of the volume of unit balls. If $\Gamma$ is finite, denote the covering space $\mathcal{K}^2\times \mathbb{R}$ or $\mathbb{S}^2\times \mathbb{R}$ as $\tilde{M}$, the covering map as $\pi$, and a fundamental domain of the cover as $M_1$. Let $B$ be any unit ball in $M^3$, we have
\begin{eqnarray*}
\operatorname{Vol}(B)&=&\operatorname{Vol}(\pi^{-1}(B)\cap M_1)
\\
&=&\int_{\pi^{-1}(B)\cap M_1}d\tilde{g}=\frac{1}{|\Gamma|}\sum_{\gamma\in\Gamma}\int_{\gamma(\pi^{-1}(B)\cap M_1)}d\tilde{g}
\\
&=&\frac{1}{|\Gamma|}\sum_{\gamma\in\Gamma}\int_{\pi^{-1}(B)\cap\gamma{ M_1}}d\tilde{g}=\frac{1}{|\Gamma|}\sum_{\gamma\in\Gamma}\int_{\gamma M_1}\chi_{\pi^{-1}(B)}d\tilde{g}
\\
&=&\frac{1}{|\Gamma|}\operatorname{Vol}(\pi^{-1}(B)),
\end{eqnarray*}
notice that we have used $\gamma(\pi^{-1}(B))=\pi^{-1}(B)$. Since $\pi^{-1}(B)$ must contain a unit ball of $\tilde{M}$, which has bounded geometry, we obtain a lower volume bound for all unit balls of $M^3$, hence bounded geometry is proved in this case.
\\

Summarizing all the cases we have discussed, we complete the proof.

\end{proof}

\textbf{Acknowledgement.} The author would like to thank his advisors, Professor Bennett Chow and Professor Lei Ni for their constant professional support. He would also like to thank Professor Natasa Sesum, whose encouragement brought his attention back to this problem.

\bibliographystyle{plain}
\bibliography{citation}

\noindent Department of Mathematics, University of California, San Diego, CA, 92093
\\ E-mail address: \verb"yoz020@ucsd.edu"

\end{document}